\documentclass[11pt]{article}
\usepackage{amsmath}
\usepackage{amssymb, amsmath}
\let\pa=\partial
\let\al=\alpha

\let\f=\frac
\let\p=\psi

\def\cR{{\cal R}}

\def\na{\nabla}

\def\p{\partial}

\def\dv{\mbox{div}}

\def\C{\mathop{\bf C\kern 0pt}\nolimits}
\def\DD{\mathop{\bf D\kern 0pt}\nolimits}
\def\K{\mathop{\bf K\kern 0pt}\nolimits}
\def\N{\mathop{\bf N\kern 0pt}\nolimits}
\def\Q{\mathop{\bf Q\kern 0pt}\nolimits}
\def\R{\mathop{\bf R\kern 0pt}\nolimits}
\def\endproof{\hphantom{MM}\hfill\llap{$\square$}\goodbreak}

\newcommand{\Z}{{\mathbf Z}}

\newcommand{\beq}{\begin{equation}}
\newcommand{\eeq}{\end{equation}}
\newcommand{\ben}{\begin{eqnarray}}
\newcommand{\een}{\end{eqnarray}}
\newcommand{\beno}{\begin{eqnarray*}}
\newcommand{\eeno}{\end{eqnarray*}}

\renewcommand{\theequation}{\thesection.\arabic{equation}}

\newtheorem{Theorem}{Theorem}[section]

\newtheorem{Proposition}[Theorem]{Proposition}
\newtheorem{Lemma}[Theorem]{Lemma}

\newtheorem{Remark}[Theorem]{Remark}

\begin{document}
\title{On the well-posedness  of the Ideal MHD equations in the Triebel-Lizorkin spaces}
\author{Qionglei Chen $^\dag$  Changxing Miao $^{\dag}$ and Zhifei Zhang
$^{\ddag}$\\[2mm]
{\small $ ^\dag$ Institute of Applied Physics and Computational Mathematics, Beijing 100088, China}\\
{\small E-mail: chen\_qionglei@iapcm.ac.cn  and   miao\_changxing@iapcm.ac.cn}\\[2mm]
{\small $ ^\ddag$ School of  Mathematical Science, Peking University, Beijing 100871, China}\\
{\small E-mail: zfzhang@math.pku.edu.cn}}

\date{}
\maketitle

\begin{abstract}
In this paper, we prove the local well-posedness for the Ideal MHD
equations in the Triebel-Lizorkin spaces and obtain blow-up
criterion of smooth solutions. Specially, we fill a gap in a step
of the proof of the local well-posedness part  for the
incompressible Euler equation in \cite{Chae1}.
\end{abstract}

 \vspace{0.2cm}


\renewcommand{\theequation}{\thesection.\arabic{equation}}
\setcounter{equation}{0}
%%%%%%%%%%%%%%%%%%%%%%%%%%%%%%%%%%%%%%%%%%%%%%
%%%%%%%%%%%%%%%%%%%%%%%%%%%%%%%%%%%%%%%%%%
\section{Introduction}
In this paper, we are concerned with  the Ideal MHD equations in
$\R^d$:
\ben\label{1.1} (\rm IMHD)\quad\left\{
\begin{aligned}
&u_t+u\cdot \nabla u=-\nabla p-\frac{1}{2}\nabla b^2+b\cdot \nabla
b,\\
&b_t+u\cdot \nabla b=b\cdot \nabla u ,\\
&\nabla\cdot u=\nabla\cdot b=0,\\
&u(0,x)=u_0(x),\quad b(0,x)=b_0(x),
\end{aligned}
\right. \een where $x\in \R^d, t\ge 0$, $u$, $b$ describes the flow
velocity vector and the magnetic field vector respectively, $p$ is a
scalar pressure, while $u_0$ and $b_0$ are
 the given initial velocity and initial magnetic field
with $\nabla\cdot u_0=\nabla\cdot b_0=0$.

Using the standard energy method \cite{Maj1}, it can be proved
that for $(u_0, b_0)\in H^s(\R^d)$, $s>\f d2+1$, there exists
$T>0$ such that the Cauchy problem (\ref{1.1}) has a unique smooth
solution $(u(t,x),b(t,x))$ on $[0,T)$ satisfying
$$(u, b)\in C([0,T);H^s)\cap C^1([0,T);H^{s-1}).$$
But whether this local solution will exist globally or lead to a
singularity in finite time is still an outstanding open problem.
Caflisch, Klapper and Steele\cite{Caf} extended
Beale-Kato-Majda criterion \cite{Bea} for the incompressible Euler equations
to the Ideal MHD equations. More precisely, they showed if the smooth
solution $(u,b)$ satisfies the following condition: \ben\label{1.2}
\int_0^T(\|\nabla\times u\|_{L^\infty}+\|\nabla\times
b\|_{L^\infty})dt<\infty, \een then the solution $(u,b)$ can be
extended beyond $t=T$, namely, for some $T<\tilde{T}, (u,b) \in
C([0,\tilde{T});H^s)\cap C^1([0,\tilde{T});H^{s-1})$. One can refer to
\cite{CCM,ZL} for the other refined criterions, and for the viscous MHD equations, some  criterions
can be found in \cite{CMZ,Wu1, Wu2, Wu3, Wu4}.

Recently, Chae studied the local well-posedness and blow-up
criterion for the incompressible Euler equations in the
Triebel-Lizorkin spaces\cite{Chae1,Chae2}. As we know,
Triebel-Lizorkin spaces are the unification of several classical
function spaces such as  Lebegue spaces $L^p(\R^d)$, Sobolev
spaces $H^s_p(\R^d)$, Lipschitz spaces $C^s(\R^d)$, and so on.
 In \cite{Chae1}, the
author first used the Littlewood-Paley operator to localize the
Euler equation to the frequency annulus $\{|\xi|\sim 2^j\}$, then
obtained an integral representation of the frequency-localized
solution on the Lagrangian coordinates by introducing a family of
particle trajectory mapping $\{X_j(\al,t)\}$ defined by
\begin{equation}\label{1.2'}
\Bigg\{\begin{aligned}
&\f {\p} {\p t}X_j(\al,t)=(S_{j-2}v)(X_j(\al,t),t)\\
&X_j(\al,0)=\al,
\end{aligned} \Bigg.
\end{equation}
where $v$ is a divergence-free velocity field and $S_{j-2}$ is a
frequency projection to the ball $\{|\xi|\lesssim 2^j\}$(see
Section 2).

 With the integral representation, one can obtain the well-posedness of the Euler equation in the
framework of the Besov spaces by standard argument, due to the
following important relation
$$\Big(\sum_{j\in\Z}2^{jsq}\|\Delta_jv(X_j(\alpha,t))\|_{L^p(\cdot d\alpha)}^q\Big)^{\frac1q}
\cong\|v\|_{\dot{B}^s_{p,q}}$$
 by  the volume-preserving property of the
mapping $\{X_j(\al,t)\}$ which is defined by \eqref{1.2'}.
However, if we work in the framework of the
Triebel-Lizorkin spaces, and  the trajectory mapping
$\{X_j(\al,t)\}$ is taken, we don't know whether the relation
\ben\label{1.2''}
\Bigl\|\Bigl(\sum_{j\in
\Z}2^{jsq}|\Delta_jv(X_j(\al,t))|^q\Bigr)^\f1q\Bigr\|_{L^p(\cdot
d\alpha)}\cong\Bigl\|\Bigl(\sum_{j\in
\Z}2^{jsq}|\Delta_jv(x)|^q\Bigr)^\f1q\Bigr\|_{L^p(\cdot
dx)}=\|v\|_{\dot F^s_{p,q}}   \een
 holds. The reason is that the mapping $\{X_j(\al,t)\}$ depends
on the index $j$, and we can't find a uniform change of the
coordinates independent of $j$ such that (\ref{1.2''}) holds. On
the other hand,   the proof of  the commutator estimate (the key
point of the proof of  the local well-posedness part)
\begin{equation}\label{a}\Big\|\Big(\sum_{j\in\Z}2^{jsq}\big|[(S_{j-2}v\cdot\na)\Delta_jv
-\Delta_j((v\cdot\na)v)](X_j(\al,t))\big|^q\Big)^{\frac1q}\Big\|_{L^p}
\le C\|\na v\|_{\infty}\|v\|_{\dot{F}^s_{p,q}}\end{equation} also
leads to some trouble  due to similar reasons.

The purpose of this paper  is to deal with  the well-posedness of
the Ideal MHD equations \eqref{1.1} in the Triebel-Lizorkin
spaces. Firstly, we  can  reduce   \eqref{1.1}   to the transport
equations  by introducing the symmetrizers.   If  we  still use
the trajectory mapping depending on $j$,  the above-mentioned
trouble will occur. In order to overcome this  difficulty, we will
introduce a different family of particle trajectory mapping
$\{X(\al,t)\}$ independent of $j$ defined by
\begin{equation}
\left\{ \begin{aligned}
&\f {\p} {\p t}X(\al,t)=v(X(\al,t),t)\\
&X(\al,0)=\al.
\end{aligned} \right.\nonumber
\end{equation}
The price to pay here is that  we have to establish  the following  commutator estimate
\begin{align}
\bigg\|\Big(\sum_{j\in\Z}2^{jsq}\big|[(v\cdot\na)\Delta_ju
-\Delta_j((v\cdot\na)u)]\big|^q\Big)^{\frac1q}\bigg\|_{L^p} \le C\big(\|\na
v\|_{L^\infty}\|u\|_{\dot F^s_{p,q}}+\|u\|_{L^\infty}\|\na v\|_{\dot
F^s_{p,q}}\big)\nonumber
\end{align}
by the paradifferential calculus, whose proof is  more
complicated since  $v$  is rougher than $S_{j-2}v$ which is the smooth low frequency cut-off of  $v$.
It is necessary to point out that
 the Maximal inequality (see Lemma 2.5 in Section
2) plays a key role  in the proof of the above inequality, which
helps us to avoid  other difficulties   arising from the change of
the coordinates.

Now we state our result as follows.

\begin{Theorem}
{\bf (i) Local-in-time Existence.}\, Let $(u_0,b_0)\in F^s_{p,q}$,
$s> \frac{d}{p}+1$, $1< p,q<\infty$ satisfying $\text{\rm div\;}u_0=\text{\rm div\;}b_0=0$.
Then there exists $T=T(\|(u_0, b_0)\|_{F^s_{p,q}})$ such that the
IMHD has a unique solution $(u,b)\in C([0,T);F^s_{p,q})$.

{\bf (ii) Blow-up Criterion.}\, The local-in-time solution $(u,b)\in
C([0,T); F^s_{p,q})$ constructed in (i) blows up at $T^*>T$ in
$F^s_{p,q}$, i.e.
$$\limsup_{t\nearrow T^*}\|(u,b)\|_{F^s_{p,q}}=+\infty,\quad T^*<\infty,$$
if and only if
\beq
\int_0^{T^*}\|(\na\times u, \na\times b)(t)\|_{\dot{F}^0_{\infty,\infty}}dt=+\infty.\label{1.3}
\eeq
\end{Theorem}
\begin{Remark}  In the case of\,  $b=0$, (IMHD) can be read  as the
incompressible Euler equations, and  what  proved in
\cite{Chae1} is a straightforward consequence of Theorem 1.1.
\end{Remark}
\begin{Remark}
Using the argument in \cite{CCM}, we can also refine blow-up criterion (\ref{1.3}) to the following form:
there exists a positive constant $M_0$ such that if
$$\lim_{\varepsilon\rightarrow0}\sup_{j\in\Z}\int_{T^*-\varepsilon}^{T^*}
\|(\Delta_j(\na\times u), \Delta_j(\na\times b))(t)\|_{L^\infty}dt\ge M_0,$$
then the smooth solution $(u,b)$ blows up at $t=T^*$.
\end{Remark}

\noindent{\bf Notation:} Throughout this paper, $C$ stands for a
``harmless" constant, and we will use the notation $A\lesssim B$ as
an equivalent to $A\le CB$, $A\approx B$ as $A\lesssim B$ and
$B\lesssim A$,  and denote $\|\cdot\|_p$ by  $L^p(\R^d)$ norm of a
function.

\setcounter{equation}{0}
\section{Preliminaries}

Let ${\cal B}=\{\xi\in\R^d,\, |\xi|\le\frac{4}{3}\}$ and ${\cal
C}=\{\xi\in\R^d,\, \frac{3}{4}\le|\xi|\le\frac{8}{3}\}$. Choose two
nonnegative smooth radial functions $\chi$, $\varphi$ supported
respectively in ${\cal B}$ and ${\cal C}$ such that
\beno
\chi(\xi)+\sum_{j\ge0}\varphi(2^{-j}\xi)=1,\quad\xi\in\R^d,\\
\sum_{j\in\Z}\varphi(2^{-j}\xi)=1,\quad\xi\in\R^d\backslash \{0\}.
\eeno We denote $\varphi_j(\xi)=\varphi(2^{-j}\xi)$, $h={\cal
F}^{-1}\varphi$ and $\tilde{h}={\cal F}^{-1}\chi$. Then the dyadic
blocks $\Delta_j$ and $S_j$ can be defined as follows
\beno
&&\Delta_jf=\varphi(2^{-j}D)f=2^{jd}\int_{\R^d}h(2^jy)f(x-y)dy, \\
&&S_jf=\sum_{k\le
j-1}\Delta_kf=\chi(2^{-j}D)f=2^{jd}\int_{\R^d}\tilde{h}(2^jy)f(x-y)dy.
\eeno
Formally, $\Delta_j=S_{ j}-S_{j-1}$  is a frequency projection
to the annulus $\{|\xi|\sim 2^j\}$, and $S_j$ is a frequency
projection to the ball $\{|\xi|\lesssim 2^j\}$. One easily verifies
that with our choice of $\varphi$
\begin{eqnarray}\label{2.1}
\Delta_j\Delta_kf\equiv0\quad i\!f\quad|j-k|\ge 2\quad and \quad
\Delta_j(S_{k-1}f\Delta_k f)\equiv0\quad i\!f\quad|j-k|\ge 5.
\end{eqnarray}

With the introduction of $\Delta_j$ and $S_j$, let us recall the definition of the Triebel-Lizorkin space. Let $s\in \R$,
$(p,q)\in[1,\infty)\times[1,\infty]$, the homogenous
Triebel-Lizorkin space $\dot {F}^s_{p,q}$ is defined by
$$\dot {F}^s_{p,q}=\{f\in {\cal Z}'(\R^d); \|f\|_{\dot
{F}^s_{p,q}}<\infty\},$$ where
$$\|f\|_{\dot{F}^s_{p,q}}=\Bigg\{\begin{array}{l}
\Big\|\big(\sum_{j\in\Z}2^{jsq}|\Delta_j f|^q\big)^{\frac 1
q}\Big\|_{p},\quad \hbox{for}\quad 1\le q<\infty,\\
\big\|\sup_{j\in \Z}(2^{js}|\Delta_jf|)\big\|_p, \quad \hbox{ for}
\quad q=\infty,
\end{array}\Bigg.
$$
and ${\cal Z}'(\R^d)$ denotes the dual space of ${\cal Z}(\R^d)=\{
f\in {\cal S}(\R^d); \partial^\alpha\hat f(0)=0; \forall \alpha\in
\N^d \,\hbox {multi-index}\}$ and can be identified by the quotient
space of ${\cal S}'/{\cal P}$ with the polynomials space ${\cal P}$.

For $s>0$, and  $(p,q)\in[1,\infty)\times[1,\infty]$, we define the
inhomogeneous Triebel-Lizorkin space  $F^s_{p,q}$ as follows
$$F^s_{p,q}=\{f\in {\cal S}'(\R^d);\, \|f\|_{
F^s_{p,q}}<\infty\},$$ where
$$\|f\|_{F^s_{p,q}}=\|f\|_{p}+\|f\|_{\dot{F}^s_{p,q}}.$$
We refer to \cite{Ber,Tri} for more details.

\begin{Lemma}(Bernstein's inequality) Let $k\in \N$. There exist a
constant $C$ independent of $f$ and  $j$ such that for all  $1\le p\le
q\le\infty$, the following inequalities hold:
\begin{align}
&{\rm supp}\hat{f}\subset\{|\xi|\lesssim 2^j\}\Rightarrow
\sup_{|\al|=k}\|\partial^\alpha f\|_q\le C2^{j{k}+j
d(\frac{1}{p}-\frac{1}{q})}\|f\|_{p},
\nonumber\\
&{\rm supp}\hat{f}\subset\{|\xi|\sim 2^j\}\Rightarrow
\|f\|_{p}\le C\sup_{|\alpha|=k}2^{-jk}\|\partial^\alpha
f\|_p.\nonumber
\end{align}
\end{Lemma}

 For the proof, see \cite{Chem2,Mey}.

\begin{Lemma}\label{lem2.2}
For any $k\in\N$, there exists a constant $C_k$ such that the
following inequality holds:
\begin{align}\label{}
C_k^{-1}\|\na^k f\|_{\dot F^s_{p,q}}\le\|f\|_{\dot F^{s+k}_{p,q}}\le
C_k\|\na^k f\|_{\dot F^s_{p,q}}.\nonumber
\end{align}
\end{Lemma}

The proof  can be found in \cite{Tri}.

\begin{Proposition}\label{Prop2.1}\cite{Chae1}
Let $s>0$, $(p,q)\in(1,\infty)\times(1,\infty]$, or $p=q=\infty$,
then there exists   a constant $C$ such that
\begin{align}
\|fg\|_{\dot F^s_{p,q}}\le C\big(\|f\|_\infty\|g\|_{\dot F^s_{p,q}}+\|g\|_{\infty}\|f\|_{\dot F^s_{p,q}}\big),\nonumber\\
\|fg\|_{F^s_{p,q}}\le
C\big(\|f\|_\infty\|g\|_{F^s_{p,q}}+\|g\|_\infty\|f\|_{F^s_{p,q}}\big).\nonumber
\end{align}
\end{Proposition}

For a locally integrable function $f$, the maximal function  $Mf(x)$ is
defined by
$$Mf(x)=\sup_{r>0}\frac{1}{|{\cal B}(x,r)|}\int_{{\cal B}(x,r)}|f(y)|dy,$$
where $|{\cal B}(x,r)|$ is the volume of the ball ${\cal B}(x,r)$
with center $x$ and radius $r$.

\begin{Lemma}\cite{Feff}\label{Lem2.3}(Vector Maximal inequality) Let $(p,q)\in (1,\infty)\times (1,\infty]$ or
$p=q=\infty$  be given. Suppose $\{f_j\}_{j\in\Z}$ is a sequence of
function in $L^p$ with the property that $\|f_j(x)\|_{\ell^q(\Z)}\in
L^p(\R^d)$. Then there holds
\begin{align}
\bigg\|\bigg(\sum_{j\in\Z}|Mf_j(x)|^q\bigg)^{\frac1q}\bigg\|_p\le C
\bigg\|\bigg(\sum_{j\in\Z}|f_j(x)|^q\bigg)^{\frac1q}\bigg\|_p.\nonumber
\end{align}
\end{Lemma}

\begin{Lemma}\label{Lem2.4}
Let $\varphi$ be an integrable function on $\R^d$, and set
$\varphi_{\varepsilon}(x)=\frac{1}{\varepsilon^d}\varphi(\frac{x}{\varepsilon})$
for $\varepsilon>0$. Suppose that the least decreasing radial
majorant of $\varphi$ is integrable; i.e. let
$$\psi(x)=\sup_{|y|\ge |x|}|\varphi(y)|,$$ and we suppose
$\int_{\R^d}\psi(x)dx=A<\infty$. Then with the same $A$, {for}
$f\in L^p(\R^d),$  $1\le p\le\infty$
$$\sup_{\varepsilon>0}|(f\ast\varphi_\varepsilon)(x)|\le AM(f)(x).$$
\end{Lemma}

The proof  can be found in  \cite{Ste}, Chap. III.

\setcounter{equation}{0}
\section{The proof of Theorem 1.1}

We divide the  proof of Theorem 1.1 into several
steps.\vspace{0.1cm}

{\bf Step 1.} A priori estimates.\vspace{0.1cm}

Let us symmetrize the equation (\ref{1.1}). Set
$$
z^+=u+b,\qquad z^-=u-b,
$$
then (\ref{1.1}) can be reduced to the system for $z^+$ and $z^-$
\ben\label{3.1}\left\{
\begin{array}{ll}
\partial_tz^++(z^-\cdot\nabla)z^+=-\nabla\pi,\\
\partial_tz^-+(z^+\cdot\nabla)z^-=-\nabla\pi,\\
\na\cdot z^+=\na\cdot z^-=0,\\
z^+(0)=z^+_0=u_0+b_0,\quad z^-(0)=z^-_0=u_0-b_0,
\end{array} \right.\quad \een
where $\pi=p+\frac{1}{2}b^2$. Taking  the operation $\Delta_k$ on both sides of (\ref{3.1}), we  get
\begin{align}\label{3.2} \Bigg\{
\begin{aligned}
\partial_t\Delta_{k}z^++z^{-}\cdot\na \Delta_{k}z^++\na\Delta_{k}\pi
&=[z^-, \Delta_k]\cdot\na z^+,\\
\partial_t\Delta_{k}z^-+z^{+}\cdot\na \Delta_{k}z^-+\na\Delta_{k}\pi
&=[z^+, \Delta_k]\cdot\na z^-,
\end{aligned}\Bigg.
\end{align}
where we  denote the commutators \beno
&[z^-, \Delta_k]\cdot\na z^+ \triangleq z^{-}\cdot\na \Delta_{k}{z^+}-\Delta_k((z^-\cdot\nabla)z^+),\\
&[z^+, \Delta_k]\cdot\na z^- \triangleq z^{+}\cdot\na
\Delta_{k}{z^-}-\Delta_k((z^+\cdot\nabla)z^-). \eeno
Let  $X^+_t(\al)$ and $X^-_t(\al)$ be the solutions of the following ordinary differential
equations:
\begin{align}\label{3.3}\left\{
\begin{aligned}
&{\pa_t}X^+_t(\al)=z^-(X^+_t(\al),t),\\
&{\pa_t}X^-_t(\al)=z^+(X^-_t(\al),t),\\
&X^+_t(\al)\big|_{t=0}=X^-_t(\al)\big|_{t=0}=\al.
\end{aligned}\right.
\end{align}
Then, it follows from \eqref{3.2} that
\begin{align}\label{3.4}\left\{
\begin{aligned}
&\frac{d}{dt}\Delta_kz^+(X^+_t(\al),t)=[z^-, \Delta_k]\cdot\na z^+(X^+_t(\al),t)-\na\Delta_{k}\pi(X^+_t(\al),t),\\
&\frac{d}{dt}\Delta_{k}z^-(X^-_t(\al),t)=[z^+, \Delta_k]\cdot\na
z^-(X^-_t(\al),t)-\na\Delta_{k}\pi(X^-_t(\al),t),
\end{aligned}\right.
\end{align}
which implies that
\begin{align}\label{3.5}
\big|\Delta_{k}z^+(X^+_t(\al),t)\big| \le &|\Delta_kz_0^+(\alpha)|+
\int_0^t\big|([z^-, \Delta_k]\cdot\na z^+)(X^+_{\tau}(\al),
\tau)\big|d\tau\nonumber\\&+ \int_0^t\big|\Delta_{k}\na
\pi(X^+_{\tau}(\al), \tau)\big|d\tau.
\end{align}
Multiplying $2^{ks}$, taking $\ell^q(\Z)$ norm   on both sides of
(\ref{3.5}), we get by using Minkowski inequality that
\begin{align}\label{3.6}
\Big(\sum_{k}|2^{ks}
\Delta_{k}z^+(X^+_{t}(\al),t)|^q\Big)^{\frac1q}\le
&\Big(\sum_{k}|2^{ks} \Delta_{k}z^+_{0}(\al)|^q\Big)^{\frac1q}+
\int_0^t\Big(\sum_{k}|2^{ks} \Delta_{k}\na
\pi(X^+_{\tau}(\al),\tau)|^q\Big)^{\frac1q}d\tau\nonumber\\&+
\int_0^t\Big(\sum_{k}|2^{ks} ([z^-, \Delta_k]\cdot\na
z^+)(X^+_{t}(\al),\tau)|^q\Big)^{\frac1q}d\tau.
\end{align}
Next, taking the $L^p$ norm with respect to $\al\in\R^d$ on both sides of \eqref{3.6},
we get by using the Minkowski inequality that
\begin{align}\label{3.7}
\Big(\int_{\R^d}\Big|&\Big(\sum_{k}|2^{ks}
\Delta_{k}z^+(X^+_t(\al),t)|^q\Big)^{\frac1q}\Big|^p
d\al\Big)^{\frac1p} \nonumber\\ \le & \|z^+_0\|_{\dot F^s_{p,q}}+
\int_0^t\Big(\int_{\R^d}\Big|\Big(\sum_{k}|2^{ks} \Delta_{k}\na \pi(X^+_{\tau}(\al),\tau)|^q\Big)^{\frac1q}\Big|^p d\al\Big)^{\frac1p}d\tau\nonumber\\
&+ \int_0^t\Big(\int_{\R^d}\Big|\Big(\sum_{k}|2^{ks} ([z^-,
\Delta_k]\cdot\na
z^+)(X^+_{\tau}(\al),\tau)|^q\Big)^{\frac1q}\Big|^p
d\al\Big)^{\frac1p}d\tau.
\end{align}
Using the fact that $X^+_{t}(\al)$ is a volume-preserving diffeomorphism due to $\dv z^+=0$,
we get from (\ref{3.7}) that
\begin{align}\label{3.8}
\|z^+(t)\|_{\dot F^s_{p,q}}&\le \|z^+_0\|_{\dot F^s_{p,q}}+
\int_0^t\|\na \pi\|_{\dot F^s_{p,q}}d\tau\nonumber\\&\quad+
\int_0^t\Big\|\big\|2^{ks}([z^-, \Delta_k]\cdot\na
z^+)\big\|_{\ell^q(k\in\Z)}\Big\|_{p}d\tau.
\end{align}
Thanks to Proposition \ref{PropA.1}, the last term on the right
side of \eqref{3.8} is dominated by
\begin{align}\label{3.9}
\int_0^t(\|\na z^+\|_\infty+\|\na z^-\|_\infty)(\|z^-\|_{\dot
F^s_{p,q}} +\|z^+\|_{\dot F^s_{p,q}})d\tau.
\end{align}
Next, we estimate the second term on the right side of \eqref{3.8}.
Taking the divergence on both sides of (\ref{3.1}), we obtain the following
representation of the pressure
\begin{align}\label{3.10}
\pi=(-\Delta)^{-1}(\pa_jz^-_i\pa_iz^+_j)=(-\Delta)^{-1}\pa_i\pa_j(z^-_iz^+_j).
\end{align}
For $l,m\in[1,d]$, we have
\begin{align}
\pa_l\pa_m\pi=(-\Delta)^{-1}\pa_l\pa_m(\pa_jz^-_i \pa_iz^+_j)={\cal
R}_l{\cal R}_m(\pa_jz^-_i \pa_iz^+_j),\nonumber
\end{align}
where $\cR_l$ denotes the Riesz transform.
Thanks to the boundedness of the Riesz transform in the homogeneous
Triebel-Lizorkin spaces \cite{FTW}, Lemma\ref{lem2.2} and Proposition \ref{Prop2.1}, we get
\begin{align}\label{3.11}
\|\na\pi\|_{\dot F^s_{p,q}}&\le
C\sum_{l,m=1}^d\|\pa_l\pa_m\pi\|_{\dot F^{s-1}_{p,q}}
\le C\|\pa_jz^-_i \pa_iz^+_j\|_{\dot F^{s-1}_{p,q}}\nonumber\\
&\le C(\|\na z^-\|_\infty\|\na z^+\|_{\dot F^{s-1}_{p,q}}+\|\na z^+\|_\infty\|\na z^-\|_{\dot F^{s-1}_{p,q}})\nonumber\\
&\le C(\|\na z^-\|_\infty\|z^+\|_{\dot F^{s}_{p,q}}+\|\na
z^+\|_\infty\|z^-\|_{\dot F^{s}_{p,q}}).
\end{align}
Plugging \eqref{3.9} and \eqref{3.11} into \eqref{3.8} yields that
\begin{align}\label{3.12}
\|z^+(t)\|_{\dot F^s_{p,q}}&\le \|z^+_0\|_{\dot F^s_{p,q}}+C
\int_0^t(\|\na z^+\|_\infty+\|\na
z^-\|_\infty)(\|z^-\|_{\dot F^s_{p,q}}
+\|z^+\|_{\dot F^s_{p,q}})d\tau.
\end{align}
Similar  argument also leads to
\begin{align}\label{3.13}
\|z^-(t)\|_{\dot F^s_{p,q}}&\le \|z^-_0\|_{\dot F^s_{p,q}}+C
\int_0^t(\|\na z^+\|_\infty+\|\na
z^-\|_\infty)(\|z^-\|_{\dot F^s_{p,q}}
+\|z^+\|_{\dot F^s_{p,q}})d\tau.
\end{align}
In order to get the inhomogeneous version of (\ref{3.12}) and (\ref{3.13}), we have to estimate the
$L^p$ norm of $(z^+, z^-)$. Multiplying the first equation of
(\ref{3.1}) by $|z^+|^{p-2}z^+$ and the second one by
$|z^-|^{p-2}z^-$, integrating the resulting equations over $\R^d$, we obtain
\begin{align}\label{3.14}
\|z^+\|_{p}+\|z^-\|_{p}\le
\|z^+_0\|_{p}+\|z^-_0\|_{p}+C\int_0^t\|\na \pi(\tau)\|_p d\tau.
\end{align}
Using \eqref{3.10} and the $L^p$-boundedness of the Riesz transform, we get
\begin{align}\label{3.15}
\|\na \pi\|_p\le C\|z^-\cdot \na z^+\|_p\le C\|\na
z^+\|_\infty\|z^-\|_p.
\end{align}
Summing up \eqref{3.12}-\eqref{3.15} yields that
\begin{align}\label{3.16}
\|z^+(t)\|_{F^s_{p,q}}&+\|z^-(t)\|_{F^s_{p,q}}\le \|z^+_0\|_{F^s_{p,q}}+\|z^-_0\|_{F^s_{p,q}}\nonumber\\
&+C\int_0^t(\|\na z^+\|_\infty+\|\na
z^-\|_\infty)(\|z^-\|_{F^s_{p,q}}
+\|z^+\|_{F^s_{p,q}})d\tau,
\end{align}
which together with the Gronwall inequality gives
\ben
\|(z^+(t),z^-(t))\|_{F^s_{p,q}}\le \|(z^+_0,z^-_0)\|_{F^s_{p,q}}\exp{\left(C\int_0^t(\|\na z^+\|_\infty+\|\na
z^-\|_\infty)d\tau\right)}.\label{3.17}
\een

{\bf Step 2.}\, Approximate solutions and uniform estimates. \vspace{0.1cm}

We construct the approximate solutions of (\ref{3.1}). Define the sequence
$\{u^{(n)}, b^{(n)}\}_{\N_0=\N\cup\{0\}}$ by solving the following systems:
\begin{align}\label{3.18}
\left\{
\begin{aligned}
&\partial_tu^{(n+1)}+u^{(n)}\cdot \nabla u^{(n+1)}-b^{(n)}\cdot
\nabla b^{(n+1)}=-\nabla
{\tilde\pi}_1^{(n+1)},\\
&\partial_tb^{(n+1)}+u^{(n)}\cdot \nabla
b^{(n+1)}-b^{(n)}\cdot\nabla u^{(n+1)}
=-\nabla{\tilde\pi}_2^{(n+1)},\\
&\nabla\cdot b^{(n+1)}=\nabla\cdot u^{(n+1)}=0,\\
&(u^{(n+1)}, b^{(n+1)})\big|_{t=0}=S_{n+2}(u_0,b_0).
\end{aligned}\right.
\end{align}

We set $(u^{(0)},b^{(0)})=(0,0)$, and
$$
{z^+}^{(n)}=u^{(n)}+b^{(n)},\qquad {z^-}^{(n)}=u^{(n)}-b^{(n)}.
$$
Then (\ref{3.18}) can be reduced to \ben{\label{3.19}}\left\{
\begin{array}{ll}
\partial_t{z^+}^{(n+1)}+({z^-}^{(n)}\cdot\nabla){z^+}^{(n+1)}=-\nabla\pi_1^{(n+1)},\\
\partial_t{z^-}^{(n+1)}+({z^+}^{(n)}\cdot\nabla){z^-}^{(n+1)}=-\nabla\pi_2^{(n+1)},\\
\na\cdot {z^+}^{(n+1)}=\na\cdot {z^-}^{(n+1)}=0, \quad \forall n\in \mathbb{N}\\
{z^+}^{(n+1)}(0)=S_{n+2}{z^+_0},\quad
{z^-}^{(n+1)}(0)=S_{n+2}{z^-_0},
\end{array} \right.\quad \een
where  $({z^+}^{(0)},{z^-}^{(0)})=(0,0)$.  Similar to the proof of  (\ref{3.16}),  we conclude that
\begin{align}\label{3.20}
&\|({z^+}^{(n+1)}(t),{z^-}^{(n+1)}(t))\|_{F^s_{p,q}}\nonumber\\
\le& \|({z^+_0},{z^-_0})\|_{F^s_{p,q}}+C\int_0^t\Big(\|(\na {z^+}^{(n)},\na{z^-}^{(n)})\|_\infty+\| (\na {z^+}^{(n+1)}, \na
{z^-}^{(n+1)}\|_\infty\Big)\nonumber\\
\qquad\quad &\times \Big(\|(\na {z^+}^{(n)},\na{z^-}^{(n)})\|_{F^s_{p,q}}+\| (\na {z^+}^{(n+1)}, \na
{z^-}^{(n+1)}\|_{F^s_{p,q}}\Big)d\tau,
\end{align}
where we used the
fact that
\begin{align}
\|(S_{n+2}{z^+_0}, S_{n+2}{z^-_0})\big\|_{F^s_{p,q}}\le\|({z^+_0},
{z^-_0})\big\|_{F^s_{p,q}}.\nonumber
\end{align}
Note that $F^{s-1}_{p,q}\hookrightarrow L^\infty$ for
$s-1>\frac{d}{p}$, (\ref{3.20}) ensures that there exists
$T_0=T_0(\|({z^+_0},{z^-_0})\big\|_{F^s_{p,q}})$ such that for any
$n$, $t\in [0,T_0]$
\begin{align}\label{3.21}
\|({z^+}^{(n)}(t),{z^-}^{(n)}(t))\|_{F^s_{p,q}}\le 2\|({z^+_0},{z^-_0})\|_{F^s_{p,q}}.
\end{align}

{\bf Step 3.} Existence.\vspace{0.1cm}

We will show that there exists a positive time $T_1(\le T_0)$
independent of $n$ such that $\{{z^+}^{(n)}, {z^-}^{(n)}\}$ is a
Cauchy sequence in $X_T^{s-1}\triangleq {\cal C}([0,T_1]; F^{s-1}_{p,q})$. For this purpose, we set
$$\delta
{z^+}^{(n+1)}={z^+}^{(n+1)}-{z^+}^{(n)},\;  \delta
{z^-}^{(n+1)}={z^-}^{(n+1)}-{z^-}^{(n)}, \;  \delta
\pi_j^{(n+1)}=\pi_j^{(n+1)}-\pi_j^{(n+1)},\; j=1, 2.$$ Using (\ref{3.19}), it is easy
to verify that the difference  $(\delta{z^+}^{(n+1)}, \delta
{z^-}^{(n+1)}, \delta \pi^{(n)})$ satisfies
\ben \left\{
\begin{aligned}\label{3.22}
&\pa_t\delta {z^+}^{(n+1)}+{z^-}^{(n)}\cdot\na\delta {z^+}^{(n+1)}=-\delta {z^-}^{(n)}\cdot\na{z^+}^{(n)}-\na\delta\pi_1^{(n+1)},\\
&\pa_t\delta {z^-}^{(n+1)}+{z^+}^{(n)}\cdot\na\delta
{z^-}^{(n+1)}=-\delta
{z^+}^{(n)}\cdot\na{z^-}^{(n)}-\na\delta\pi_2^{(n+1)} ,\\&(\delta
{z^+}^{(n+1)},  \delta
{z^-}^{(n+1)})\big|_{t=0}=\Delta_{n+1}(z^+_0,z^-_0).
\end{aligned}\right.
\een
Applying $\Delta_k$ to the first equation of \eqref{3.22}, we get
\begin{align}\label{3.23}
\pa_t\Delta_k\delta {z^+}^{(n+1)}+{z^-}^{(n)}\cdot\na\Delta_k\delta
{z^+}^{(n+1)}=& [{z^-}^{(n)},
\Delta_k]\cdot\na\delta{z^+}^{(n+1)}\nonumber\\& -\Delta_k(\delta
{z^-}^{(n)}\cdot\na{z^+}^{(n)})-\na\Delta_k\delta\pi_1^{(n+1)}.
\end{align}
Exactly as in the proof of (\ref{3.8}), we get
\begin{align}\label{3.24}
\|\delta{z^+}^{(n+1)}\|_{\dot F^{s-1}_{p,q}}\le&
C\|\Delta_{n+1}z^+_0\|_{\dot F^{s-1}_{p,q}}
+\int_0^t\Big\|\big\|2^{k(s-1)}([{z^-}^{(n)}, \Delta_k]\cdot\na\delta {z^+}^{(n+1)})(\al,\tau)\big\|_{\ell^q(\Z)}\Big\|_{p}d\tau\nonumber\\
&+\int_0^t\|\delta {z^-}^{(n)}\cdot\na{z^+}^{(n)}(\tau)\|_{\dot
F^{s-1}_{p,q}}d\tau+ \int_0^t\|\na\delta\pi_1^{(n+1)}(\tau)\|_{\dot
F^{s-1}_{p,q}}d\tau.
\end{align}
Thanks to the Fourier support of $\Delta_{n+1}z^+_0$, we have
\begin{align}\label{3.25}
\|\Delta_{n+1}z^+_0\|_{\dot F^{s-1}_{p,q}}\le
C2^{-(n+1)}\|z^+_0\|_{\dot F^s_{p,q}}.
\end{align}
Using Proposition 4.1 and the embedding $F^{s-1}_{p,q}\hookrightarrow L^\infty$,
the second term on the right side of
\eqref{3.24} is dominated by
\begin{align}\label{3.26}
&\|\na{z^-}^{(n)}\|_{\infty}\|\delta{z^+}^{(n+1)}\|_{\dot
F^{s-1}_{p,q}}+\|\delta{z^+}^{(n+1)}\|_\infty \|\na
{z^-}^{(n)}\|_{\dot F^{s-1}_{p,q}} \nonumber\\ \le&
C\|{z^-}^{(n)}\|_{F^s_{p,q}}\|\delta{z^+}^{(n+1)}\|_{F^{s-1}_{p,q}}.
\end{align}
Thanks to Proposition 2.1, the third term on the right hand side
of \eqref{3.24} is dominated by
\begin{align}\label{3.27}
&\|\delta {z^-}^{(n)}\|_\infty\|\na {z^+}^{(n)}\|_{\dot
F^{s-1}_{p,q}}+\|\delta {z^-}^{(n)}\|_{\dot
F^{s-1}_{p,q}}\|\na{z^+}^{(n)}\|_\infty\nonumber\\
\le&
C\|\delta{z^-}^{(n)}\|_{F^{s-1}_{p,q}}\|{z^+}^{(n)}\|_{F^s_{p,q}}.
\end{align}
Taking the divergence on both sides of \eqref{3.22}, we get
$$
\delta\pi_1^{(n+1)}=\pa_j(-\Delta)^{-1}(\delta {z^-_i}^{(n)}
\pa_i {z^+_j}^{(n)})+\pa_i(-\Delta)^{-1}(\pa_j{z^-_i}^{(n)}\delta
{z^+_j}^{(n+1)}).
$$
Hence, we have
\begin{align}\pa_l\delta\pi_1^{(n+1)}={\cal R}_l{\cal R}_j(\delta
{z^-_i}^{(n)}\pa_i {z^+_j}^{(n)})+{\cal R}_l{\cal
R}_i(\pa_j{z^-_i}^{(n)}\delta {z^+_j}^{(n+1)}),\nonumber
\end{align}
which together with Proposition \ref{Prop2.1} and the boundedness of the Riesz
transform in the homogeneous
Triebel-Lizorkin spaces gives
\begin{align}\label{3.28}
\|\na\delta\pi_1^{(n+1)}\|_{\dot
F^{s-1}_{p,q}}\lesssim&  \|\delta {z^-_i}^{(n)}\pa_i
{z^+_j}^{(n)}\|_{\dot F^{s-1}_{p,q}}
+\|\pa_j{z^-_i}^{(n)}\delta {z^+_j}^{(n+1)}\|_{\dot F^{s-1}_{p,q}}\nonumber\\
\lesssim &\|\delta {z^-}^{(n)}\|_\infty\|\na {z^+}^{(n)}\|_{\dot
F^{s-1}_{p,q}}+
\|\delta {z^-}^{(n)}\|_{\dot F^{s-1}_{p,q}}\|\na {z^+}^{(n)}\|_\infty\nonumber\\
&+\|\na {z^-}^{(n)}\|_\infty\|\delta {z^+}^{(n+1)}\|_{\dot
F^{s-1}_{p,q}} +\|\na {z^-}^{(n)}\|_{\dot F^{s-1}_{p,q}}\|\delta
{z^+}^{(n+1)}\|_\infty\nonumber\\ \lesssim& \|\delta
{z^-}^{(n)}\|_{F^{s-1}_{p,q}}\|{z^+}^{(n)}\|_{F^s_{p,q}}+\|{z^-}^{(n)}\|_{F^s_{p,q}}\|\delta
{z^+}^{(n+1)}\|_{F^{s-1}_{p,q}}.
\end{align}
By summing up \eqref{3.24}-\eqref{3.28}, we get
\begin{align}\label{3.29}
\|\delta{z^+}^{(n+1)}\|_{\dot F^{s-1}_{p,q}}\le
C2^{-(n+1)}\|z^+_0\|_{\dot F^s_{p,q}}+C\int_0^t
\Big(&\|{z^-}^{(n)}\|_{F^s_{p,q}}\|\delta{z^+}^{(n+1)}\|_{F^{s-1}_{p,q}}\nonumber\\
&+\|\delta
{z^-}^{(n)}\|_{F^{s-1}_{p,q}}\|{z^+}^{(n)}\|_{F^s_{p,q}}\Big)dt.
\end{align}
Now, we estimate the $L^p$ norm of $\delta{z^+}^{(n+1)}$. Multiplying $|\delta
{z^+}^{(n+1)}|^{p-2}\delta {z^+}^{(n+1)}$ on both sides of the first
equation of \eqref{3.22}, and integrating the resulting equations over $\R^d$, we obtain
\begin{align}
\|\delta {z^+}^{(n+1)}(t)\|_p\le &\|\Delta_{n+1}z^+_0\|_p+
\int_0^t\|\delta {z^-}^{(n)}\cdot\na{z^+}^{(n)}(\tau)\|_pd\tau+
\int_0^t\|\na \delta\pi_1^{(n+1)}(\tau)\|_pd\tau\nonumber\\ \le &
2^{-(n+1)}\|z^+_0\|_{\dot F^s_{p,q}}+C\int_0^t\|\delta
{z^-}^{(n)}\|_p\|\na{z^+}^{(n)}(\tau)\|_\infty d\tau\nonumber\\&+
C\int_0^t\|\na {z^-}^{(n)}\|_\infty\|\delta {z^+}^{(n+1)}\|_pd\tau,\nonumber
\end{align}
which together with \eqref{3.29} gives
\begin{align}\label{3.30}
\|\delta{z^+}^{(n+1)}\|_{F^{s-1}_{p,q}}\le
C2^{-(n+1)}\|z^+_0\|_{F^s_{p,q}}+C\int_0^t
\Big(&\|{z^-}^{(n)}\|_{F^s_{p,q}}\|\delta{z^+}^{(n+1)}\|_{F^{s-1}_{p,q}}\nonumber\\
&+\|\delta
{z^-}^{(n)}\|_{F^{s-1}_{p,q}}\|{z^+}^{(n)}\|_{F^s_{p,q}}\Big)dt.
\end{align}
Exactly as in the proof of \eqref{3.30}, we also have
\begin{align}\label{3.31}
\|\delta{z^-}^{(n+1)}\|_{F^{s-1}_{p,q}}\le
C2^{-(n+1)}\|z^-_0\|_{F^s_{p,q}}+C\int_0^t
\Big(&\|{z^+}^{(n)}\|_{F^s_{p,q}}\|\delta{z^-}^{(n+1)}\|_{F^{s-1}_{p,q}}\nonumber\\
&+\|\delta
{z^+}^{(n)}\|_{F^{s-1}_{p,q}}\|{z^-}^{(n)}\|_{F^s_{p,q}}\Big)dt.
\end{align}
Adding  up \eqref{3.30} and \eqref{3.31}, we obtain
\begin{align}
&\big\|\big(\delta{z^+}^{(n+1)},
\,\delta{z^-}^{(n+1)}\big)\big\|_{F^{s-1}_{p,q}} \nonumber\\ &\quad\lesssim
2^{-(n+1)}(\|z^+_0\|_{F^s_{p,q}}+\|z^-_0\|_{F^s_{p,q}})\nonumber\\&\qquad+ T\sup_{t\in[0,T]}
\big\|\big({z^+}^{(n)},\,{z^-}^{(n)}\big)\big\|_{F^s_{p,q}}\big\|\big(\delta{z^+}^{(n+1)},\,\delta{z^-}^{(n+1)}\big)\big\|_{F^{s-1}_{p,q}}
\nonumber\\&\qquad+T\sup_{t\in[0,T]}\big\|\big({z^+}^{(n)},\,{z^-}^{(n)}\big)\big\|_{F^s_{p,q}}
\big\|\big(\delta {z^+}^{(n)},\,\delta
{z^-}^{(n)}\big)\big\|_{F^{s-1}_{p,q}},\nonumber
\end{align}
which together with \eqref{3.21} yields that
\begin{align}\label{x}
\big\|\big(\delta{z^+}^{(n+1)},
\,\delta{z^-}^{(n+1)}\big)\big\|_{X^{s-1}_T}\le& C_12^{-(n+1)}
+C_1T\big\|\big(\delta{z^+}^{(n+1)},
\,\delta{z^-}^{(n+1)}\big)\big\|_{X^{s-1}_T}\nonumber\\&+
C_1T\big\|\big(\delta{z^+}^{(n)},
\,\delta{z^-}^{(n)}\big)\big\|_{X^{s-1}_T},
\end{align}
where $C_1=C_1(\|(z^+_0,\,z^-_0)\|_{F^s_{p,q}})$.  Thus, if $C_1T\le\frac14$, then
\begin{align}
\big\|\big(\delta{z^+}^{(n+1)},
\,\delta{z^-}^{(n+1)}\big)\big\|_{X^{s-1}_T}\le C_12^{-n}
+2C_1T\big\|\big(\delta{z^+}^{(n)},
\,\delta{z^-}^{(n)}\big)\big\|_{X^{s-1}_T}.\nonumber
\end{align}
This implies that
\begin{align}
\big\|\big(\delta{z^+}^{(n+1)},
\,\delta{z^-}^{(n+1)}\big)\big\|_{X^{s-1}_T}\le&
2C_12^{-(n+1)}.\nonumber
\end{align}
Thus, $\{{z^+}^{(n)}, {z^-}^{(n)}\}_{n\in\N_0}$ is a
Cauchy sequence in $X_{T_1}^{s-1}$. By the
standard argument, for $T_1\le\min\{T_0,\frac1{4C_1}\}$, the limit $(z^+,z^-)\in X_{T_1}^{s}$ solves the equation
(\ref{3.1}) with the initial data $(z^+_0,\,z^-_0)$. Moreover,
$(z^+,z^-)$ satisfies
\begin{align}\label{}
&\|(z^+,z^-)(t)\|_{L^\infty_{T_1}(F^s_{p,q})}\le
C\|(z^+_0,\,z^-_0)\|_{F^s_{p,q}},\nonumber
\end{align}
which implies $(u, b)$  is  a solution of \eqref{1.1} with the
initial data $(u_0,\,b_0)\in F^s_{p,q}$, and
\begin{align}\label{}
&\|(u,\,b)(t)\|_{L^\infty_{T_1}(F^s_{p,q})}\le
C\|(u_0,\,b_0)\|_{F^s_{p,q}}.\nonumber
\end{align}

{\bf The proof of the uniqueness.} Consider $({z^+}', {z^-}')\in
C_{T_1}(F^s_{p,q})$ is another solution to (\ref{3.1}) with the same
initial data. Let $\delta z^+=z^+-{z^+}'$ and $\delta
z^-=z^--{z^-}'$. Then $(\delta z^+,\delta z^-)$ satisfies the
following equations \ben \nonumber\left\{
\begin{aligned}
&\pa_t\delta z^++(z^-\cdot\na)\delta z^+=-(\delta z^-\cdot\na)z^+-\na(\pi-\pi'),\\
&\pa_t\delta z^-+(z^+\cdot\na)\delta z^-=-(\delta z^+\cdot\na)z^--\na(\pi-\pi'),\\
&\nabla\cdot \delta z^+=\nabla\cdot\delta z^-=0.
\end{aligned}\right.
\een In the same way as deriving in \eqref{x}, we obtain
\begin{align}
\big\|\big(\delta{z^+}, \,\delta{z^-}\big)\big\|_{X^{s-1}_T}\le&
C_2T\big\|\big(\delta{z^+},
\,\delta{z^-}\big)\big\|_{X^{s-1}_T}\nonumber
\end{align}
for sufficiently small $T$.  This implies that $(\delta z^+,
\,\delta z^-)\equiv0$, i.e.,  $(z^+, \,z^-)\equiv ({z^+}', \,{z^-}')$.
\vspace{.2cm}

{\bf Blow-up Criterion.}  By means of  Proposition 1.1 in
\cite{Chae1} and
$$\|(\na z^+,\na z^-)\|_{\dot
F^0_{\infty,\infty}}\lesssim \|(\na\times z^+,\na\times
z^-)\|_{\dot F^0_{\infty,\infty}},$$
we have
$$\|(\na z^+,\na z^-)\|_\infty\lesssim \Big(1+\|(\na z^+,\na
z^-)\|_{\dot F^0_{\infty,\infty}}\big(\log\big(1+\|(\na\times z^+,\na\times
z^-)\|_{F^{s-1}_{p,q}}\big)+1\big)\Big)$$
Plugging the above estimates into \eqref{3.16} then by Gronwall's lemma yields that
$$\|(z^+,z^-)\|_{F^s_{p,q}}\le \|(z^+_0,
z^-_0)\|_{F^s_{p,q}}\exp\Big[C \exp[C\int_0^t(1+\|(\na\times
z^+,\na\times z^-)\|_{\dot F^0_{\infty,\infty}}d\tau]\Big]$$ which
implies the blow-up criterion. This finishes the proof of the
Theorem 1.1.

\vspace*{4mm}

\section{Appendix}

Let us recall the para-differential calculus which enables us to
define a generalized product between distributions, which is
continuous in many functional spaces where the usual product does
not make sense (see \cite{Bon}). The para-product between $u$ and
$v$ is defined by$$T_uv\triangleq\sum_{j\in\Z}S_{j-1}u\Delta_jv.$$
We then have the following formal decomposition: \beq\label{A.1}
uv=T_uv+T_vu+R(u,v), \eeq with
$$R(u,v)=\sum_{j\in\Z}\Delta_ju\widetilde{\Delta}_jv\quad\mbox{and}\quad
\widetilde{\Delta}_j=\Delta_{j-1}+\Delta_j+\Delta_{j+1}.$$
 The
decomposition (\ref{A.1}) is called the Bony's para-product
decomposition.
\begin{Proposition}\label{PropA.1}
Let $(p,q)\in(1,\infty)\times(1,\infty]$, or $p=q=\infty$, and $f$
be a solenoidal vector field. Then for $s>0$
\begin{align}\label{A.2}
\Big\|\big\|2^{ks}([f, \Delta_k]\cdot\na
g)\big\|_{\ell^q(\Z)}\Big\|_{p} \lesssim \big(\|\na
f\|_\infty\|g\|_{\dot F^s_{p,q}}+\|\na g\|_\infty\|f\|_{\dot
F^s_{p,q}}\big).
\end{align}
or for $s>-1$
\begin{align}\label{A.3}
\Big\|\big\|2^{ks}([f, \Delta_k]\cdot\na
g)\big\|_{\ell^q(\Z)}\Big\|_{p} \lesssim \big(\|\na
f\|_\infty\|g\|_{\dot F^s_{p,q}}+\|g\|_\infty\|\na f\|_{\dot
F^s_{p,q}}\big).
\end{align}
\end{Proposition}

\noindent{\it Proof.} By the Einstein convention on the summation
over repeated indices $i\in[1,d]$, and the Bony's paraproduct
decomposition we decompose \begin{align} [f, \Delta_k]\cdot\na g= [f_i,
\Delta_k]\pa_i g&= [T_{f_i}, \Delta_k]\pa_i g+T'_{\Delta_k\pa_i
g}{f}_i -\Delta_k(T_{\pa_i g}f_i)-\Delta_k(R(f_i,\pa_i
g))\nonumber\\&\triangleq I+II+III+IV\nonumber,\end{align} where $T'_uv$ stands for
$T_uv+R(u,v)$. Thank to  the support condition (\ref{2.1}), we
rewrite
\begin{align}\label{A.4}
|I|&=\Big|\sum_{k'\sim k}[S_{k'-1}f_i, \Delta_k]\pa_i
\Delta_{k'}g\Big|\nonumber\\&=\Big|\sum_{k'\sim k}\int_{\R^d}
\Big(S_{k'-1}f_{i}(x)-S_{k'-1}f_{i}(y)\Big)2^{kd}h(2^k(x-y))\pa_i
\Delta_{k'}g(y) dy\Big|,
\end{align}
where $k'\sim k$ stands for $|k'-k|\le4$. Integrate by part and use
$\dv f=0$, the integrand in (\ref{A.4}) is
\begin{align}
\Big(S_{k'-1}f_i(x)-S_{k'-1}f_i(y)\Big)2^{k(d+1)}(\pa_i h)(2^k(x-y))
\Delta_{k'} g(y)\nonumber
\end{align}
which was dominated by
\begin{align}\label{A.5}
\|\na S_{k'-1}f\|_\infty 2^k|x-y|2^{kd}|\na
h(2^k(x-y))||\Delta_{k'}g(y)|.
\end{align}
Recall $h(x)\in{\cal S}(\R^d)$, it is easy to see that  $|x\na
h(x)|$ satisfies Lemma $\ref{Lem2.4}$, so (\ref{A.5}) is less than
\begin{align}\label{A.6}
C\|\na S_{k'-1}f\|_\infty M(|\Delta_{k'}g(\cdot)|)(x).
\end{align}
Multiplying $2^{ks}$ on both sides of (\ref{A.4}), taking
$\ell^q(\Z)$ norm  then taking $L^p$ norm and putting (\ref{A.6})
into the resulting inequality, we have
\begin{align}\label{A.7}
\big\|\|2^{ks}|I(x)|\|_{\ell^q(\Z)}\big\|_{p}&\lesssim \|\na
S_{k'-1}f\|_{\infty}\Big\| \big\|\sum_{k'\sim
k}2^{(k-k')s}M(2^{k's}|\Delta_{k'}g(\cdot)|)(x)\big\|_{\ell^q(\Z)}\Big\|_p\nonumber\\&\lesssim
\|\na
f\|_\infty\Big\|\big\|M(2^{k's}|\Delta_{k'}g(\cdot)|)(x)\big\|_{\ell^q(\Z)}\Big\|_p\nonumber\\&\lesssim
\|\na
f\|_\infty\big\|\big\|2^{ks}|\Delta_{k}g(x)|\big\|_{\ell^q(\Z)}\big\|_p\lesssim
\|\na f\|_\infty\|g\|_{\dot F^s_{p,q}},
\end{align}
where  we used  Lemma {\ref{Lem2.3} in the third inequality.
Let us turn to the term $II$, thanks to the definition of $II$,
\begin{align}\label{A.8}
|II|=\Big|\sum_{k'\ge k-2}S_{k'+2}\pa_i \Delta_{k}g
\Delta_{k'}f_i(x)\Big| \le\sum_{k'\ge k-2}\|\na
\Delta_{k}g\|_\infty|\Delta_{k'}f(x)|.
\end{align}
Then thanks to the convolution inequality for series,  we get for
$s>0$,
\begin{align}\label{A.9}
\big\|\|2^{ks}|II(x)|\|_{\ell^q(\Z)}\big\|_{p}&\lesssim\|\na
\Delta_{k}g\|_\infty\Big\| \big\|\sum_{k'\ge
k-2}2^{(k-k')s}2^{k's}|\Delta_{k'}f(x)|\big\|_{\ell^q(\Z)}\Big\|_p\nonumber\\&\lesssim
\|\na
\Delta_{k}g\|_\infty\Big\|\big\|2^{-ks}\chi_{\{k\ge-2\}}\big\|_{\ell^1(\Z)}\big\|2^{k's}|\Delta_{k'}f(x)|\big\|_{\ell^q(\Z)}\Big\|_p
\nonumber\\&\lesssim \|\na
g\|_\infty\big\|\big\|2^{ks}|\Delta_{k}f(x)|\big\|_{\ell^q(\Z)}\big\|_p
\lesssim \|\na g\|_\infty\|f\|_{\dot F^s_{p,q}}.
\end{align}
For the term $III$,
\begin{align}\label{A.10}
|III|&=\Big|\sum_{k'\sim k}\Delta_k(S_{k'-1}{\pa_i
g}\Delta_{k'}f_i)\Big|\lesssim \sum_{k'\sim k}\big|M(S_{k'-1}{\pa_i
g}\Delta_{k'}f_i)(x)\big|\nonumber\\& \lesssim\sum_{k'\sim
k}\big|M(|\Delta_{k'}f|)(x)\big|\|S_{k'-1}\na g\|_\infty.
\end{align}
Using (\ref{A.10}) and in the same way as leading to (\ref{A.7})
yields
\begin{align}\label{A.11}
\big\|\|2^{ks}|III(x)|\|_{\ell^q(\Z)}\big\|_{p}&\lesssim \|\na
S_{k'-1}g\|_{\infty}\Big\| \big\|\sum_{k'\sim
k}2^{(k-k')s}M(2^{k's}|\Delta_{k'}f|)(x)\big\|_{\ell^q(\Z)}\Big\|_p\nonumber\\&\lesssim
\|\na
g\|_\infty\Big\|\big\|M(2^{k's}|\Delta_{k'}f|)(x)\big\|_{\ell^q(\Z)}\Big\|_p\nonumber\\&\lesssim
\|\na
g\|_\infty\big\|\big\|2^{ks}|\Delta_{k}f(x)|\big\|_{\ell^q(\Z)}\big\|_p\lesssim
\|\na g\|_\infty\|f\|_{\dot F^s_{p,q}}.
\end{align}
In view of $\dv f=0$ and integrating by part, we have
\begin{align}\label{A.12}
|IV|=&\Big|\sum_{k'\ge k-3}\Delta_k(\Delta_{k'}f_i\pa_i
\widetilde{\Delta}_{k'}g)\Big|
=\Big|\sum_{k'\ge k-3}\int_{\R^d}2^{kd}h(2^k(x-y))\Delta_{k'}f_i(y)\pa_i \widetilde{\Delta}_{k'}g(y)dy\Big|\nonumber\\
=&\Big|\sum_{k'\ge k-3}\int_{\R^d}2^{kd+k}(\pa_i h)(2^k(x-y))\Delta_{k'}f_i(y)\widetilde{\Delta}_{k'}g(y)dy\Big|\nonumber\\
\lesssim&\sum_{k'\ge
k-3}2^kM(\Delta_{k'}f\widetilde{\Delta}_{k'}g)(x)
\lesssim\sum_{k'\ge
k-3}2^kM(\widetilde{\Delta}_{k'}g)(x)\|\Delta_{k'}f\|_\infty.
\end{align}
The convolution inequality for series and Lemma \ref{Lem2.3} allow us
to give that for $s+1>0$,
\begin{align}\label{A.13}
\big\|\|2^{ks}|IV(x)|\|_{\ell^q(\Z)}\big\|_{p}&\lesssim \|\na
\Delta_{k'}f\|_{\infty}\Big\| \big\|\sum_{k'\ge
k-3}2^{(k-k')(s+1)}M(2^{k's}\widetilde{\Delta}_{k'}g)(x)\big\|_{\ell^q(\Z)}\Big\|_p\nonumber\\&\lesssim
\|\na
f\|_\infty\Big\|\big\|M(2^{k's}\widetilde{\Delta}_{k'}g)(x)\big\|_{\ell^q(\Z)}\Big\|_p\nonumber\\&\lesssim
\|\na
f\|_\infty\big\|\big\|2^{ks}|\widetilde{\Delta}_{k}g(x)|\big\|_{\ell^q(\Z)}\big\|_p\lesssim
\|\na f\|_\infty\|g\|_{\dot F^s_{p,q}}.
\end{align}
Summing up \eqref{A.7}, \eqref{A.9}, \eqref{A.11} and \eqref{A.13},
we get the desired inequality \eqref{A.2}.

In order to prove  the inequality  \eqref{A.3},  we only indicate how to get the bound on  $II$ and
$III$ since  $I$ and $IV$ can be treated as above.   We  estimate the term $II$ as
\begin{align}
|II|=\Big|\sum_{k'\ge k-2}S_{k'+2}\pa_i \Delta_{k}g
\Delta_{k'}f_i(x)\Big| \le\sum_{k'\ge
k-2}2^k\|\Delta_{k}g\|_\infty|\Delta_{k'}f(x)|.\nonumber
\end{align}
Then thanks to the convolution inequality for series,  we get for
$s+1>0$,
\begin{align}
\big\|\|2^{ks}|II(x)|\|_{\ell^q(\Z)}\big\|_{p}&\lesssim\|\Delta_{k}g\|_\infty\Big\|
\big\|\sum_{k'\ge
k-2}2^{(k-k')(s+1)}2^{k'(s+1)}|\Delta_{k'}f(x)|\big\|_{\ell^q(\Z)}\Big\|_p\nonumber\\&\lesssim
\|g\|_\infty\Big\|\big\|2^{-k(s+1)}\chi_{\{k\ge-2\}}\big\|_{\ell^1(\Z)}\big\|2^{k'(s+1)}|\Delta_{k'}f(x)|\big\|_{\ell^q(\Z)}\Big\|_p
\nonumber\\&\lesssim \|g\|_\infty\|f\|_{\dot
F^{s+1}_{p,q}}.\nonumber
\end{align}
Let's turn to the term $III$,
\begin{align}
|III|\lesssim\sum_{k'\sim
k}\big|M(|\Delta_{k'}f|)(x)\big|2^{k'}\|S_{k'-1}
g\|_\infty.\nonumber
\end{align}
Arguing similarly as in deriving \eqref{A.11} yields that
\begin{align}
\big\|\|2^{ks}|III(x)|\|_{\ell^q(\Z)}\big\|_{p}&\lesssim \|
S_{k'-1}g\|_{\infty}\Big\| \big\|\sum_{k'\sim
k}2^{(k-k')s}M(2^{k'(s+1)}|\Delta_{k'}f|)(x)\big\|_{\ell^q(\Z)}\Big\|_p\nonumber\\&\lesssim
\| g\|_\infty\|f\|_{\dot F^{s+1}_{p,q}}.\nonumber
\end{align}
Thus the desired inequality \eqref{A.3} is obtained.  \endproof

\vspace*{4mm}

\textbf{Acknowledgements}  The authors  would like  to thank
Professor D. Chae for his helpful discussion and comments. The
authors are deeply grateful to the referee for his
 invaluable comments and suggestions which helped improve the paper greatly.
  Q. Chen, C. Miao and Z.Zhang  were
supported by the NSF of China
under grant No.10701012, No.10725102 and No.10601002.

\end{document}